\documentclass{amsart}

\usepackage{amsmath,amssymb,mathtools}
\usepackage{microtype}
\usepackage[hidelinks]{hyperref}

\newtheorem{theorem}{Theorem}[section]
\newtheorem{proposition}[theorem]{Proposition}
\newtheorem{corollary}[theorem]{Corollary}
\newtheorem{lemma}[theorem]{Lemma}
\theoremstyle{remark}
\newtheorem{remark}[theorem]{Remark}

\usepackage{xcolor}
\newcommand{\Ecolor}[1]{{\color{black}{#1}}}
\newcommand{\N}{\mathbb{N}}
\newcommand{\E}{\mathcal{E}}
\newcommand{\D}{\mathcal{D}}
\newcommand{\F}{\mathcal{F}}
\newcommand{\T}{\mathcal{T}}
\newcommand{\Res}{\operatorname*{Res}}

\title[On the Density of Repeated Entries in Pascal's Triangle]{On the Density of Repeated Entries in Pascal's Triangle}
\author{Amir Leshem}
\author{Eitan Sayag}
\thanks{Amir Leshem was partially supported by ISF grant no. 2197/22. Eitan Sayag was partially supported by ISF grant no. 1781/23 and BSF grant no. 2022240.}

\begin{document}
\begin{abstract}
Let $N(a)$ denote the number of occurrences of the integer $a>1$ in
Pascal's triangle, and set
\[
 \E_H=\{a\geq3:N(a)>2H\},
 \qquad
 A_H(x)=\#(\E_H\cap[1,x]).
\]
For every fixed $K\geq2$ we prove the expansion
\[
 A_1(x)=\sum_{k=2}^{K}(k!x)^{1/k}
 +O_K\!\left(x^{1/(K+1)}\right).
\]
Consequently, the Dirichlet series
\[
 \D_1(s)=\sum_{a\in\E_1}a^{-s}
\]
extends meromorphically to $\Re s>0$, with simple poles precisely at
$s=1/k$, $k\geq2$, of residue $(k!)^{1/k}/k$.

\Ecolor{
For every $H\geq2$, $x\ge 3$ and every $\varepsilon>0$, we prove 
\[
A_H(x)\le A_2(x)\leq\exp\left(\frac{2\log x}{(\log\log x)^{1/4}}\right)\ll_{\varepsilon} x^{\varepsilon}.
\]
}

Hence 
\Ecolor{
\[
 \D_H(s)=\sum_{a\in\E_H}a^{-s}
\]
}
is holomorphic for $\Re s>0$ when $H\geq2$, and
\[
 \lambda(\E_1)=\frac12,
 \qquad
 \lambda(\E_H)=0\quad(H\geq2).
\]
\Ecolor{Here, for a set $S\subseteq\N$, its exponent of convergence is
\[
\lambda(S)=\inf\left\{\sigma>0:
 \sum_{n\in S}n^{-\sigma}<\infty\right\}.
\]
}
We also determine the logarithmic-order counting law for the classical
Lind--Singmaster family of coefficients of multiplicity at least six.
\end{abstract}

\maketitle

\section{Introduction and historical context}
Singmaster \cite{Singmaster1971} conjectured that there \Ecolor{is} an absolute bound on the occurrences of any number \Ecolor{in} Pascal's triangle. The largest presently known multiplicity is eight, attained by
\[
 3003=\binom{3003}{1}=\binom{78}{2}=\binom{15}{5}=\binom{14}{6}
\]
and the four symmetric occurrences.  On the other hand, Lind and Singmaster found an infinite Fibonacci family of values occurring at least six times \cite{Lind1968,Singmaster1975}. Formally, 
for an integer $a>1$, let
\begin{equation}\label{eq:N-def}
N(a)=\#\left\{(n,k)\in\N^2:1\leq k<n,\ \binom{n}{k}=a\right\}.
\end{equation}
Every $a\geq3$ has the two edge representations
\[
 a=\binom{a}{1}=\binom{a}{a-1},
\]
so $N(a)\geq2$.  

The first general upper bound, $N(a)=O(\log a)$, is due to Singmaster.  Abbott, Erd\H{o}s and Hanson improved it to
\[
 N(a)=O\left(\frac{\log a}{\log\log a}\right)
\]
and proved that the average and normal order of $N(a)$ is $2$ \cite{AEH}.  Kane subsequently obtained the strongest known unconditional global estimate
\begin{equation}\label{eq:Kane-bound}
 N(a)=O\left(\frac{(\log a)(\log\log\log a)}{(\log\log a)^3}\right)
\end{equation}
\cite{Kane2004,Kane2007}.  More recently, Matom\"aki, Radziwi\l\l, Shao, Tao and Ter\"av\"ainen proved Singmaster's conjecture in the interior region
\[
 \exp\bigl((\log n)^{2/3+\delta}\bigr)
 \leq k\leq
 n-\exp\bigl((\log n)^{2/3+\delta}\bigr)
\]
for every fixed $\delta>0$: a sufficiently large value has at most four occurrences in this region, or at most two in either half of Pascal's triangle \cite{MRSTT}.

There is also a substantial literature on collisions between two prescribed diagonals,
\begin{equation}\label{eq:fixed-diagonal-equation}
 \binom{u}{r}=\binom{v}{s},
 \qquad 2\leq r<s.
\end{equation}
Avanesov determined the solutions for $(r,s)=(2,3)$.  The case $(2,4)$ was solved independently by de Weger and Pint\'er, the case $(3,4)$ by de Weger, and the case $(2,5)$ by Bugeaud, Mignotte, Siksek, Stoll and Tengely \cite{Avanesov,deWeger1996,Pinter1995,deWeger1997,BMSST}.  Kiss proved finiteness for $(r,s)=(2,p)$ with $p$ an odd prime \cite{Kiss1988}.  In full generality, Beukers, Shorey and Tijdeman proved that \eqref{eq:fixed-diagonal-equation} has only finitely many positive integral solutions for every fixed pair $r<s$ \cite{BST}.  Their argument ultimately invokes Siegel's theorem and is ineffective.

In contrast, the present \Ecolor{paper} studies a different, distributional aspect of the problem.  For an integer $H\geq1$ and a real number $x\geq1$, define
\begin{equation}\label{eq:E-A-def}
 \E_H=\{a\geq3:N(a)>2H\},
 \qquad
 A_H(x)=\#(\E_H\cap[1,x]).
\end{equation}
For complex $s$ with sufficiently large real part, define
\begin{equation}\label{eq:D-H-def}
 \D_H(s)=\sum_{a\in\E_H}a^{-s}.
\end{equation}
For a set $S\subseteq\N$, its exponent of convergence is
\begin{equation}\label{eq:lambda-def}
\lambda(S)=\inf\left\{\sigma>0:
 \sum_{n\in S}n^{-\sigma}<\infty\right\}.
\end{equation}
Throughout, $\varepsilon$ denotes an arbitrary positive real number.  We write $X\ll_{\theta}Y$ if $|X|\leq C_{\theta}Y$ for a positive constant depending only on the parameter or collection of parameters $\theta$.  Thus $X\ll_{H,\varepsilon}Y$ permits dependence on $H$ and $\varepsilon$, but not on $x$.

The first exceptional set has previously been studied by Abbott et al. \cite{AEH} who proved that $A_1(x)=O\left(\sqrt{x}\right)$, and later using ideal
convergence by Gubo et al. \cite{GuboEtAl2003}.  For $0<q\leq1$, let
\[
 \mathcal I_c^{(q)}
 =\left\{S\subseteq\N:\sum_{a\in S}a^{-q}<\infty\right\}.
\]
Gubo et al. \cite{GuboEtAl2003} proved that
\[
 \mathcal I_c^{(q)}\text{-}\lim_{a\to\infty}N(a)=2
 \quad\Longleftrightarrow\quad q>\frac12,
\]
or equivalently,
\[
 \sum_{a\in\E_1}a^{-q}<\infty
 \quad\Longleftrightarrow\quad q>\frac12.
\]
  Thus their result already determines
$\lambda(\E_1)=1/2$, including divergence at the critical exponent. Subsequently, 
T\'oth, Bukor, Filip and Zsilinszky \cite{TothEtAl2021} placed this conclusion
in the finer scale of exponent-of-convergence ideals
\[
 \mathcal I_{<q}=\{S\subseteq\N:\lambda(S)<q\},
 \qquad
 \mathcal I_{\leq q}=\{S\subseteq\N:\lambda(S)\leq q\},
\]
and proved
\[
 \mathcal I_{\leq 1/2}\text{-}\lim_{a\to\infty}N(a)=2.
\]
The ratio-block framework summarized and developed by Miska and T\'oth
\cite{MiskaToth2023} further relates exponent-of-convergence ideals to
asymptotic ratio distributions and regular variation.  We return to this
connection below.

None of these works studies the higher exceptional sets $\E_H$ for
$H>1$ or the higher-order expansion of $A_1(x)$.

 In this paper we first analyze the higher exceptional sets $\E_H$, $H\geq2$, where membership forces collisions between distinct nontrivial diagonals.  The strongest of these estimates, specialized to $H=2$, then shows uniformly that all pairwise and higher-order intersections among any fixed collection of diagonals are negligible. Using these results, we improve prior analysis of the first exceptional set $\E_1$.

\begin{theorem}[The higher exceptional sets]\label{thm:EH-main}
For every integer $H\geq2$, the following two statements hold.
\begin{enumerate}
\item The Bombieri--Pila determinant method gives, for every $\varepsilon>0$,
\begin{equation}\label{eq:main-BP}
 A_H(x)\ll_{H,\varepsilon}
 x^{1/(H(H+1))+\varepsilon},
\end{equation}
with an effective implied constant.  Consequently,
\begin{equation}\label{eq:main-BP-series}
 \D_H(s)\ \text{converges absolutely for}\ 
 \Re s>\frac1{H(H+1)},
\end{equation}
and hence is holomorphic in this half-plane, with an effective proof of convergence.
\item \Ecolor{For every $x\geq3$,\begin{equation}{\label{eq:main-E2-subpoly}}
 A_H(x)\le A_2(x)\leq\exp\left(\frac{2\log x}{(\log\log x)^{1/4}}\right) \ll_{\varepsilon} x^{\varepsilon}.
\end{equation}}

 Consequently, $\D_H(s)$ converges absolutely and locally uniformly, and hence is holomorphic, throughout the half-plane $\Re s>0$.  Equivalently, $\lambda(\E_H)=0$.
\end{enumerate}
\end{theorem}

The two parts of Theorem~\ref{thm:EH-main} use, respectively, absolute irreducibility together with Bombieri--Pila and Remond's bound on the number of rational points for every fixed pair of diagonals.  The weaker elementary diagonal-growth estimate is separated from the main theorem and proved as Lemma~\ref{lem:effective}.  It isolates the basic mechanism behind all three estimates: sufficiently many distinct nontrivial lower indices force a representation on a sufficiently high diagonal.  This gives the clearest intuition for the later collision arguments and is also used in the analysis of $A_1(x)$ and $\D_1(s)$, in particular to establish the initial half-plane of absolute convergence and justify the Abel-summation passage to an improper integral.

\begin{theorem}[The first exceptional set]\label{thm:E1-main}
The set $\E_1$ satisfies the following statements.
\begin{enumerate}
\item For every fixed integer $K\geq2$,
\begin{equation}\label{eq:main-E1}
 A_1(x)
 =\sum_{k=2}^{K}(k!x)^{1/k}
 +O_K\!\left(x^{1/(K+1)}\right).
\end{equation}
In particular,
\[
 A_1(x)=\sqrt{2x}+(6x)^{1/3}+O(x^{1/4}),
\]
and
\begin{equation}\label{eq:main-E1-lambda}
 \lambda(\E_1)=\frac12.
\end{equation}
\item The function $\D_1$, initially defined for $\Re s>1/2$, admits a meromorphic continuation to $\Re s>0$.  Its poles there are precisely the simple real poles
\[
 s=\frac1k,\qquad k=2,3,4,\ldots,
\]
and
\begin{equation}\label{eq:main-D1-residues}
 \Res_{s=1/k}\D_1(s)=\frac{(k!)^{1/k}}{k}.
\end{equation}
There are no nonreal poles.
\end{enumerate}
\end{theorem}

Theorem~\ref{thm:E1-main} uses the second estimate in Theorem~\ref{thm:EH-main} with $H=2$.  It implies that the total overcount from all pairwise and higher-order intersections among the first $K-1$ nontrivial diagonals is $O_{K,\varepsilon}(x^\varepsilon)$.  Thus the fixed diagonals contribute independently to every finite-order expansion of $A_1(x)$: the $k$th diagonal supplies the term $(k!x)^{1/k}$ and, analytically, the simple pole $s=1/k$ of $\D_1(s)$.  The resulting dichotomy in the open right half-plane is sharp: every $\D_H$, $H\geq2$, is holomorphic, whereas $\D_1$ is meromorphic with the simple real poles $1/k$.

\section{Combinatorial effective estimate of $A_H(x)$}

For $a\geq3$, let
\begin{equation}\label{eq:M-def}
 M(a)=\#\left\{(n,k)\in\N^2:
 2\leq k\leq\frac n2,\ \binom{n}{k}=a
 \right\}.
\end{equation}
The edge representations are not included in $M(a)$.  Every noncentral left-half representation contributes two symmetric occurrences to $N(a)$, while a central representation contributes one.  Therefore
\begin{equation}\label{eq:N-M}
 2+M(a)\leq N(a)\leq2+2M(a).
\end{equation}
For fixed $k$, the function $n\mapsto\binom{n}{k}$ is strictly increasing for $n\geq k$.  Thus the lower indices counted by $M(a)$ are distinct.

We use the elementary inequalities
\begin{equation}\label{eq:binomial-lower}
 \binom{n}{k}
 =\prod_{j=0}^{k-1}\frac{n-j}{k-j}
 \geq\left(\frac nk\right)^k,
 \qquad 1\leq k\leq n,
\end{equation}
and, for $n\geq2k$,
\begin{equation}\label{eq:central-lower}
 \binom{n}{k}\geq\binom{2k}{k}\geq2^k.
\end{equation}

\begin{lemma}\label{lem:pair-count}
For integers $K\geq2$ and real $x\geq2$, let
\[
 P_K(x)=\#\left\{(n,k)\in\N^2:
 K\leq k\leq\frac n2,\ \binom{n}{k}\leq x
 \right\}.
\]
Then
\begin{equation}\label{eq:pair-sum}
 P_K(x)\leq
 \sum_{k=K}^{\lfloor\log_2x\rfloor}kx^{1/k}.
\end{equation}
In particular,
\begin{equation}\label{eq:pair-uniform}
 P_K(x)\ll x^{1/K}(\log x)^2
\end{equation}
uniformly in $K$, and, for each fixed $K$,
\begin{equation}\label{eq:pair-fixed}
 P_K(x)\ll_Kx^{1/K}.
\end{equation}
All implied constants are effective.
\end{lemma}

\begin{proof}
If $\binom{n}{k}\leq x$, then \eqref{eq:binomial-lower} gives $n\leq kx^{1/k}$, while \eqref{eq:central-lower} gives $k\leq\log_2x$.  This proves \eqref{eq:pair-sum}.  Bounding each $x^{1/k}$ by $x^{1/K}$ gives \eqref{eq:pair-uniform}.  If $K$ is fixed, separate the term $k=K$ and use
\[
 \sum_{k=K+1}^{\lfloor\log_2x\rfloor}kx^{1/k}
 \ll x^{1/(K+1)}(\log x)^2
 =o_K(x^{1/K}).
\]
After enlarging the constant for bounded $x$, this proves \eqref{eq:pair-fixed}.
\end{proof}
\begin{lemma}[Elementary diagonal growth]\label{lem:effective}
For every fixed integer $H\geq1$ we have \Ecolor{$A_H(x) \leq P_{H+1}(x).$ Hence,}
\begin{equation}\label{eq:effective-count}
 A_H(x)\ll_Hx^{1/(H+1)}.
\end{equation}
Consequently,
\begin{equation}\label{eq:effective-series}
 \sum_{a\in\E_H}a^{-s}<\infty
 \qquad\left(\Ecolor{\Re(s)}>\frac1{H+1}\right).
\end{equation}
\end{lemma}

\begin{proof}
If $a\in\E_H$, then $N(a)>2H$.  By the upper bound in \eqref{eq:N-M}, $M(a)\geq H$.  Hence $a$ has $H$ distinct nontrivial left-half lower indices, all at least $2$, and the largest is at least $H+1$.  Every $a\in\E_H\cap[1,x]$ is therefore represented by a pair counted by $P_{H+1}(x)$. 
Thus $A_H(x) \leq P_{H+1}(x).$ Equation \eqref{eq:effective-count}  follows from
Equation \eqref{eq:pair-fixed} of  
Lemma~\ref{lem:pair-count}.

For \eqref{eq:effective-series}, divide the positive integers into dyadic intervals.  If $\Ecolor{\gamma:=\Re(s)}>1/(H+1)$, then
\Ecolor{
\begin{align*}
 |\sum_{a\in\E_H}a^{-s}| \leq
  \sum_{a\in\E_H}a^{-\gamma}
 \leq\sum_{j\geq1}2^{-j \gamma}A_H(2^{j+1})
 \ll_H\sum_{j\geq1}2^{-j(\gamma-1/(H+1))}<\infty.
\end{align*}
}
\end{proof}

\section{Geometric bounds on $A_H(x)$ for $H \geq 2$}\label{sec:collisions}
In this section we improve the combinatorial bounds, using techniques for bounding the number of integer solutions on algebraic curves.

For fixed integers $2\leq r<s$, define the bounded-height collision count
\begin{equation}\label{eq:I-rs}
 I_{r,s}(x)=\#\left\{(u,v)\in\N^2:
 u\geq2r,\ v\geq2s,\
 \binom{u}{r}=\binom{v}{s}\leq x
 \right\}.
\end{equation}
After clearing denominators, the collision equation is the affine plane curve
\begin{equation}\label{eq:collision-curve}
 F_{r,s}(U,V)=s!(U)_r-r!(V)_s=0,
\end{equation}
where $(Z)_m=Z(Z-1)\cdots(Z-m+1)$.  This curve has degree $s$.

\begin{lemma}\label{lem:absolute-irreducibility}
For every pair of integers $2\leq r<s$, the polynomial
$F_{r,s}(U,V)$ in \eqref{eq:collision-curve} is absolutely irreducible.
\end{lemma}

\begin{proof}
Put 
\[
 \Ecolor{Q_m}(X)=X(X+1)\cdots(X+m-1).
\]
After the affine changes of variables
\[
 X=U-r+1,
 \qquad
 Y=V-s+1,
\]
we have $(U)_r=\Ecolor{Q}_r(X)$ and $(V)_s=\Ecolor{Q}_s(Y)$.  Hence
\[
 F_{r,s}(U,V)
 =s!\left(\Ecolor{Q}_r(X)-\frac{r!}{s!}\Ecolor{Q}_s(Y)\right).
\]
By \cite[Theorem~2.1]{BST}, the polynomial
\[
 \Ecolor{Q}_m(X)-\lambda \Ecolor{Q}_n(Y),
 \qquad m\leq n,\quad \lambda\in\mathbb C^\times,
\]
is reducible over $\mathbb C$ only in the cases
\[
 m=n,\ \lambda=1;
 \qquad
 m=n\ \text{odd},\ \lambda=-1;
 \qquad
 (m,n,\lambda)=\left(2,4,\frac14\right).
\]
We apply this with $m=r$, $n=s$, and $\lambda=r!/s!$.  Since $r<s$, the first two cases are impossible.  In the remaining exceptional pair $(r,s)=(2,4)$, one has
\[
 \frac{r!}{s!}=\frac1{12}\ne\frac14.
\]
Thus no exceptional case occurs, and $F_{r,s}$ is absolutely irreducible.
\end{proof}

\begin{proposition}[Bombieri--Pila bound for two diagonals]\label{prop:BP-pair}
For fixed integers $2\leq r<s$ and every $\varepsilon>0$,
\begin{equation}\label{eq:BP-pair}
 I_{r,s}(x)\ll_{r,s,\varepsilon}
 x^{1/(rs)+\varepsilon}.
\end{equation}
The implied constant is effective.
\end{proposition}

\begin{proof}
If $\binom{u}{r}\leq x$, then \eqref{eq:binomial-lower} gives
\[
 u\leq r x^{1/r}.
\]
Likewise,
\[
 v\leq s x^{1/s}\leq s x^{1/r}.
\]
Thus all relevant integral points of \eqref{eq:collision-curve} lie in a square of side
\(
 B\ll_{r,s}x^{1/r}.
\)
The Bombieri--Pila theorem states that an absolutely irreducible plane curve of degree $d\geq2$ has
\(
 O_{d,\eta}(B^{1/d+\eta})
\)
integral points in such a square \cite{BombieriPila1989}.  Applying it with $d=s$ and $\eta=r\varepsilon$ proves \eqref{eq:BP-pair}.
\end{proof}

The pairwise estimate can be summed because, for a fixed multiplicity threshold $H$, only finitely many small pairs are needed; the remaining large diagonals are handled by Lemma~\ref{lem:pair-count}.

\begin{proposition}\label{prop:BP-global}
For every integer $H\geq2$ and every $\varepsilon>0$,
\begin{equation}\label{eq:BP-global}
 A_H(x)\ll_{H,\varepsilon}
 x^{1/(H(H+1))+\varepsilon}.
\end{equation}
The implied constant is effective.
\end{proposition}

\begin{proof}
Let $a\in\E_H$.  As in the proof of Lemma~\ref{lem:effective}, $a$ has (at least) $H$ distinct nontrivial left-half lower indices, all at least $2$, and the largest is at least $H+1$. 
In particular, the largest, $k_H$, is at least $H+1$.

Set $D=H(H+1)-1$.  If $k_H>D$, then $a$ has a representation counted by $P_{D+1}(x)$, and Lemma~\ref{lem:pair-count} gives
\[
 P_{D+1}(x)\ll_H x^{1/(H(H+1))}.
\]
Otherwise the two largest indices form a pair $H\leq r<s\leq D$, so $a$ is counted by $I_{r,s}(x)$.  Consequently,
\begin{equation}\label{eq:BP-sum-pairs}
 A_H(x)
 \leq
 \sum_{H\leq r<s\leq D}I_{r,s}(x)
 +P_{D+1}(x).
\end{equation}
The sum contains only $O_H(1)$ pairs.  By Proposition~\ref{prop:BP-pair}, its largest exponent occurs at $(r,s)=(H,H+1)$.  Hence
\[
 \sum_{H\leq r<s\leq D}I_{r,s}(x)
 \ll_{H,\varepsilon}
 x^{1/(H(H+1))+\varepsilon},
\]
which, together with the tail estimate, proves \eqref{eq:BP-global}.
\end{proof}

\begin{corollary}\label{cor:BP-series}
For every integer $H\geq2$,
\begin{equation}\label{eq:BP-series}
 \sum_{a\in\E_H}a^{-s}<\infty
 \qquad\left(\Ecolor{\Re(s)}>\frac1{H(H+1)}\right),
\end{equation}
with an effective proof.
\end{corollary}

\begin{proof}
Fix $\Ecolor{\gamma:=\Re(s)}>1/(H(H+1))$ and choose $\varepsilon>0$ such that
\[
 \frac1{H(H+1)}+\varepsilon<\Ecolor{\gamma}.
\]
Dyadic decomposition and Proposition~\ref{prop:BP-global} give a convergent geometric series.
\end{proof}

For fixed integers $2\leq r<s$, set
\begin{equation}\label{eq:J-rs}
 J(r,s)=\#\left\{(u,v)\in\N^2:
 u\geq2r,\ v\geq2s,\ \binom{u}{r}=\binom{v}{s}
 \right\}.
\end{equation}
Because $u\mapsto\binom{u}{r}$ and $v\mapsto\binom{v}{s}$ are strictly increasing, $J(r,s)$ is also the number of common values of the two left-half diagonals.  Beukers, Shorey and Tijdeman \cite{BST} proved that $J(r,s)<\infty$, by an ineffective argument.  We use instead R\'emond's explicit bound for the number of rational points on a plane curve of genus at least two \cite{Remond2010}, together with the genus computation of \cite{BST}.

\begin{lemma}\label{lem:E2-M}
For every $a\geq3$,
\begin{equation}\label{eq:E2-M}
 a\in\E_2
 \quad\Longleftrightarrow\quad
 M(a)\geq2.
\end{equation}
\end{lemma}

\begin{proof}
If $M(a)\leq1$, then \eqref{eq:N-M} gives $N(a)\leq4$, so $a\notin\E_2$.  Conversely, assume $M(a)\geq2$.  Choose two left-half representations with distinct lower indices.  At most one is central, since $\binom{2k}{k}$ is strictly increasing in $k$.  They therefore contribute at least three non-edge occurrences, in addition to the two edge occurrences.  Thus $N(a)\geq5$, and $a\in\E_2$.
\end{proof}
\color{black}
\begin{lemma}\label{lem:genus}
Let $2\leq r<s$ with $(r,s)\notin\{(2,3),(2,4)\}$.  The smooth projective model of the curve \eqref{eq:collision-curve} has genus at least $2$.
\end{lemma}

\begin{proof}
With $X=U-r+1$ and $Y=V-s+1$ the curve is $X(X+1)\cdots(X+r-1)=\lambda\,Y(Y+1)\cdots(Y+s-1)$ with $\lambda=r!/s!$, and it is irreducible by Lemma~\ref{lem:absolute-irreducibility}.  By \cite[Theorem~2.2]{BST} its genus exceeds one unless $(r,s)$ is one of $(2,5)$, $(2,6)$, $(2,8)$, $(3,4)$ and $\lambda$ takes one of the special values listed there.  For these pairs $\lambda=1/60$, $1/360$, $1/20160$, $1/4$, and none of them is on the list.
\end{proof}

\begin{proposition}\label{prop:ineffective}
For all $x\geq3$,
\begin{equation}\label{eq:E2-subpoly}
 A_2(x)\leq\exp\left(\frac{2\log x}{(\log\log x)^{1/4}}\right).
\end{equation}
Consequently, for every $H\geq2$ and every $\varepsilon>0$,
\begin{equation}\label{eq:EH-subpoly}
 A_H(x)\ll_{\varepsilon}x^\varepsilon,
\end{equation}
uniformly in $H$, with an effective implied constant.
\end{proposition}

\begin{proof}
Write $\ell=\log\log x$.  We treat three cases. If $\ell\leq16$, then $\ell^{1/4}\leq2$, and \eqref{eq:E2-subpoly} follows from $A_2(x)\leq x$.

Let $16\leq\ell\leq1024$.  From Lemma~\ref{lem:effective} we have $A_2(x)\leq P_3(x)$, and \eqref{eq:pair-sum} gives $P_3(x)\leq x^{1/3}(\log_2x)^2$.  Since $\log x\geq e^{16}$ we have $(\log_2x)^2\leq x^{1/50}$.
Indeed, writing \(u=\log x\ge e^{16}\), we have \(u\log 2\ge 20000\).
Thus, by \(e^t\ge t^2/2\) for \(t\ge0\),
\[
x^{1/100}=e^{u/100}
\ge \frac{u^2}{20000}
\ge \frac{u}{\log 2}
= \log_2 x.
\]
Squaring gives \((\log_2 x)^2\le x^{1/50}\).

Now $\frac13+\frac1{50}<2/1024^{1/4}\leq2/\ell^{1/4}$ and
we obtain:
$$
A_2(x) \leq P_3(x) \leq x^{\frac{1}{3}} \log_2(x)^2
\leq x^{\frac{1}{3}+\frac{1}{50}}
\leq x^{2\ell^{-1/4}}
$$
This gives \eqref{eq:E2-subpoly}.

\Ecolor{Suppose now that $\ell\geq1024$. Fix $2\leq r<s$ with
$(r,s)\notin\{(2,3),(2,4)\}$, and write $F=F_{r,s}$. By
Lemmas~\ref{lem:absolute-irreducibility} and~\ref{lem:genus}, $F$ is absolutely
irreducible of degree $s$ and its curve has genus at least two.
Its integer coefficients have absolute value at most $M=r!s!$, so
$h_\infty(F)\leq\log M$. For $K=\mathbb Q$, the height factor in
\cite[Proposition~10.1]{Remond2010} is
$H_2(F)=\max\{h_\infty(F),1\}\log(h_\infty(F)+2)$.
Since the pairs counted by $J(r,s)$ are rational zeros of $F$, R\'emond's
bound gives
\[
 J(r,s)\leq\exp\!\bigl(5^{s^4-64}(\log M)\log(\log M+2)\bigr).
\]
Here $\log M<2\log(s!)\leq2s\log s<2s^2$ and
$\log(\log M+2)\leq2\log(2s)\leq2s$. As $4\leq5^{64}$, it follows that
$J(r,s)\leq\exp(s^3 5^{s^4})$.}
For $(r,s)=(2,3)$ and $(2,4)$ the equation has been solved completely: $J(2,3)=3$ by \cite{Avanesov} and $J(2,4)=1$ by \cite[Theorem~1]{deWeger1996}; see also \cite{Pinter1995}. Hence, for every $D\geq4$,
\begin{equation}\label{eq:C-D}
 C_D:=\sum_{2\leq r<s\leq D}J(r,s)\leq4+D^2\exp\bigl(D^3\,5^{D^4}\bigr).
\end{equation}
If $a\in\E_2\cap[1,x]$, Lemma~\ref{lem:E2-M} supplies two distinct nontrivial left-half lower indices.  If both are at most $D$, then $a$ is counted by one of the fixed pairs contributing to $C_D$.  Otherwise $a$ has a representation with lower index at least $D+1$.  By \eqref{eq:pair-sum},
\begin{equation}\label{eq:E2-split}
 A_2(x)\leq C_D+P_{D+1}(x)\leq C_D+x^{1/(D+1)}(\log_2x)^2 .
\end{equation}

Write \(L=\log x=e^\ell\) and \(T=L/\ell^{1/4}\), so that the
right-hand side of (30) is \(e^{2T}\). Take
\[
D=\left\lfloor\left(\frac{\ell}{2\log 5}\right)^{1/4}\right\rfloor.
\]
Since \(\ell\ge 1024\), we have \(D\ge 4\). The two inequalities
defining the floor give
\[
5^{D^4}\le e^{\ell/2}=\sqrt L,
\qquad
\frac{L}{D+1}<(2\log 5)^{1/4}T<\frac32 T.
\]
Moreover, \(D^3\le D^4\le\ell\). Since
\(\sqrt L=e^{\ell/2}\ge \ell^2/8\ge \ell^{5/4}\), it follows that
\[
D^3 5^{D^4}\le \ell\sqrt L\le T.
\]
Substituting these bounds into (32) and (33), we obtain
\[
A_2(x)
\le 4+D^2e^T+(\log_2 x)^2e^{3T/2}
\le \bigl(4+D^2+(\log_2 x)^2\bigr)e^{3T/2}.
\]
Now \(4+D^2\le 4+\ell\le L^2\) and \(\log_2 x=L/\log 2\le 2L\),
so the factor in parentheses is at most \(5L^2\le L^3=e^{3\ell}\).
Finally, \(T\ge\ell\sqrt L\ge 6\ell\), whence
\[
A_2(x)\le e^{3\ell+3T/2}\le e^{2T}.
\]
This proves (30).

Finally, the right-hand side of \eqref{eq:E2-subpoly} is $\ll_\varepsilon x^\varepsilon$, and $\E_H\subseteq\E_2$ for every $H\geq2$.  This gives \eqref{eq:EH-subpoly} uniformly in $H$.
\end{proof}

\color{black}
\begin{corollary}\label{cor:all-s}
For every integer $H\geq2$ and every real $s>0$,
\begin{equation}\label{eq:all-s}
 \sum_{a\in\E_H}a^{-s}<\infty.
\end{equation}
Equivalently, $\lambda(\E_H)=0$ for every $H\geq2$.
\end{corollary}

\begin{proof}
Choose $0<\varepsilon<s$.  By Proposition~\ref{prop:ineffective} and dyadic decomposition,
\[
 \sum_{a\in\E_H}a^{-s}
 \leq\sum_{j\geq1}2^{-js}A_H(2^{j+1})
 \ll_{\Ecolor{s,\varepsilon}}
 \sum_{j\geq1}2^{-j(s-\varepsilon)}<\infty.
\]

\end{proof}

\begin{corollary}\label{cor:DH-holomorphic}
For every integer $H\geq2$, the function $\D_H$ converges absolutely and
locally uniformly, and hence is holomorphic, in the half-plane $\Re s>0$.
\end{corollary}
\begin{proof}
Let $C$ be a compact subset of the half-plane $\Re s>0$, and put
\[
 \sigma_0=\min_{s\in C}\Re s>0.
\]
For $s\in C$ and $a\in\E_H$,
\[
 |a^{-s}|=a^{-\Re s}\leq a^{-\sigma_0}.
\]
The majorant is summable by Corollary~\ref{cor:all-s}.  The Weierstrass
$M$-test therefore gives uniform convergence on $C$, and holomorphy follows.
\end{proof}
\color{black}
We now reformulate our results in the language of
Miska and T\'oth \cite{MiskaToth2023}. 
For an infinite set $A=\{a_1<a_2<\cdots\}\subseteq\mathbb N$, define its
$n$th ratio block by
\[
 \mathcal R_n(A)
 =
 \left(
 \frac{a_1}{a_n},\frac{a_2}{a_n},\ldots,\frac{a_n}{a_n}
 \right).
\]
Given a distribution function $g$ on $[0,1]$, write $A\in U(g)$ if
\[
 \frac{1}{n}
 \#\left\{
 1\leq j\leq n:\frac{a_j}{a_n}<x
 \right\}
 \longrightarrow g(x)
\]
at every continuity point $x$ of $g$.  Write $A\in I(g)$ if $A$ is
contained in some set belonging to $U(g)$.  We shall use
\[
 c_0(x)=
 \begin{cases}
 0,&x=0,\\
 1,&0<x\leq1.
 \end{cases}
\]
Thus $A\in U(c_0)$ means that, asymptotically, almost all elements of each
ratio block are arbitrarily close to $0$.

Let $H\geq2$.  Since $\lambda(\E_H)=0$, the results of
Miska and T\'oth \cite{MiskaToth2023} imply that
\[
 \E_H\in I(c_0).
\]

If $\E_H$ is infinite and its ratio-block sequence has a unique
asymptotic distribution function, then that function must be $c_0$.
\section{The asymptotic formula for $A_1(x)$}
For every integer $k\geq2$, put
\begin{equation}\label{eq:Tk-def}
 \T_k=\left\{\binom{n}{k}:n\geq2k\right\}.
\end{equation}
Thus $\T_k$ is the set of values on the left half of the $k$th
nontrivial diagonal.  Every element of $\T_k$ has its two edge
occurrences and at least one nontrivial occurrence, and hence belongs
to $\E_1$.  Conversely, every element of $\E_1$ has a nontrivial
left-half representation.  Therefore
\begin{equation}\label{eq:E1-union}
 \E_1=\bigcup_{k\geq2}\T_k.
\end{equation}

\begin{proposition}\label{prop:E1-asymptotic}
For every fixed integer $K\geq2$,
\begin{equation}\label{eq:E1-expansion}
 A_1(x)
 =
 \sum_{k=2}^{K}(k!x)^{1/k}
 +O_K\!\left(x^{1/(K+1)}\right).
\end{equation}
In particular,
\begin{equation}\label{eq:E1-asymptotic}
 A_1(x)=\sqrt{2x}+(6x)^{1/3}+O(x^{1/4}),
\end{equation}
and $\lambda(\E_1)=1/2$.
\end{proposition}

\begin{proof}
Fix $K\geq2$.  For each $2\leq k\leq K$, write
\[
 X_k=(k!x)^{1/k}.
\]
The inequalities
\[
 \frac{(n-k+1)^k}{k!}
 \leq \binom{n}{k}
 \leq \frac{n^k}{k!}
\]
show that $n\leq X_k$ implies $\binom{n}{k}\leq x$, whereas
$\binom{n}{k}\leq x$ implies $n\leq X_k+k-1$.  Since
$n\mapsto\binom{n}{k}$ is strictly increasing, it follows that
\begin{equation}\label{eq:Tk-count}
 \#(\T_k\cap[1,x])=(k!x)^{1/k}+O_k(1).
\end{equation}

For $a\leq x$, define
\[
 m_K(a)=\#\{k\in\{2,\ldots,K\}:a\in\T_k\}.
\]
Thus
\[
 \sum_{k=2}^{K}\#(\T_k\cap[1,x])
 =\sum_{a\leq x}m_K(a).
\]
The complete overcount relative to the union of the first $K-1$
nontrivial diagonals is
\begin{equation}\label{eq:Delta-K}
 \Delta_K(x)
 =\sum_{a\leq x}\bigl(m_K(a)-1\bigr)_+.
\end{equation}
\Ecolor{Here we denote, for $z \in \mathbb{Z}$, by $z_+$ the value of $z$ if it is positive and zero otherwise. Thus,  $$z_+:=max(z,0).$$
}
If $m_K(a)\geq2$, then $a$ has at least two distinct nontrivial
left-half representations.  By Lemma~\ref{lem:E2-M}, this 
implies $a\in\E_2$.  Since $m_K(a)-1\leq K-2$,
Proposition~\ref{prop:ineffective} gives, for every $\varepsilon>0$,
\begin{equation}\label{eq:Delta-K-bound}
 \Delta_K(x)
 \leq (K-2)A_2(x)
 \ll_{K,\varepsilon}x^\varepsilon.
\end{equation}
In particular, \eqref{eq:Delta-K-bound} controls simultaneously all
pairwise and higher-order intersections among
$\T_2,\ldots,\T_K$.

Every element of $\E_1$ not lying on one of the first $K-1$
nontrivial diagonals has a left-half representation with lower index
at least $K+1$.  Hence Lemma~\ref{lem:pair-count} gives
\begin{equation}\label{eq:E1-tail}
 \#\left(\left(\E_1\setminus\bigcup_{k=2}^{K}\T_k\right)
 \cap[1,x]\right)
 \leq P_{K+1}(x)
 \ll_K x^{1/(K+1)}.
\end{equation}
Since
\[
 \#\left(\left(\bigcup_{k=2}^{K}\T_k\right)\cap[1,x]\right)
 =\sum_{k=2}^{K}\#(\T_k\cap[1,x])-\Delta_K(x),
\]
we obtain
\[
 \left|A_1(x)-\sum_{k=2}^{K}\#(\T_k\cap[1,x])\right|
 \leq P_{K+1}(x)+\Delta_K(x).
\]
Choose $0<\varepsilon<1/(K+1)$.  Equations
\eqref{eq:Tk-count}, \eqref{eq:Delta-K-bound}, and
\eqref{eq:E1-tail} now yield \eqref{eq:E1-expansion}.  Taking $K=3$
gives \eqref{eq:E1-asymptotic}.

Finally, Lemma~\ref{lem:effective} gives convergence of
$\sum_{a\in\E_1}a^{-s}$ for $s>1/2$.  Conversely,
$\T_2\subseteq\E_1$, and
\[
 \sum_{a\in\E_1}a^{-s}
 \geq\sum_{n\geq4}\binom{n}{2}^{-s},
\]
whose convergence threshold is $s=1/2$.  Hence
$\lambda(\E_1)=1/2$.
\end{proof}

The complete asymptotic expansion in
Proposition~\ref{prop:E1-asymptotic} determines all singularities of
$\D_1$ in the right half-plane.

\begin{theorem}\label{thm:D1-continuation}
The Dirichlet series $\D_1(s)$, initially absolutely convergent for
$\Re s>1/2$, admits a meromorphic continuation to the half-plane
$\Re s>0$.  Its poles there are precisely the simple real poles
\begin{equation}\label{eq:D1-poles}
 s=\frac1k,\qquad k=2,3,4,\ldots.
\end{equation}
\Ecolor{The residue of $\D_1(s)$ at these poles is given by:}
\begin{equation}\label{eq:D1-residue}
 \Res_{s=1/k}\D_1(s)=\frac{(k!)^{1/k}}{k}
 \qquad(k\geq2).
\end{equation}
In particular, $\D_1$ admits no
meromorphic continuation to any neighborhood of the origin. 
\end{theorem}

\begin{proof}
We first record explicitly the Abel summation formula used below.
For every $X\geq1$ and every complex number $s$,
\begin{equation}\label{eq:D1-Abel-finite}
 \sum_{\substack{a\in\E_1\\ a\leq X}}a^{-s}
 =A_1(X)X^{-s}
 +s\int_1^X A_1(x)x^{-s-1}\,dx.
\end{equation}
Indeed, this is the usual partial summation identity applied to the
counting function $A_1$.  Since
$A_1(X)\ll X^{1/2}$ by Lemma~\ref{lem:effective}, the boundary
term $A_1(X)X^{-s}$ tends to zero as $X\to\infty$ whenever
$\Re s>1/2$.  Hence, in the half-plane of absolute convergence,
\begin{equation}\label{eq:D1-Abel}
 \D_1(s)=s\int_1^\infty A_1(x)x^{-s-1}\,dx.
\end{equation}

Fix an integer $K\geq2$.  By
Proposition~\ref{prop:E1-asymptotic},
\[
 A_1(x)=\sum_{k=2}^{K}(k!)^{1/k}x^{1/k}+R_K(x),
 \qquad
 R_K(x)=O_K\!\left(x^{1/(K+1)}\right).
\]
Substituting this expansion into \eqref{eq:D1-Abel} yields
\begin{equation}\label{eq:D1-finite-expansion}
 \D_1(s)
 =
 \sum_{k=2}^{K}(k!)^{1/k}\frac{s}{s-1/k}+H_K(s),
\end{equation}
where
\begin{equation}\label{eq:D1-HK}
 H_K(s)=s\int_1^\infty R_K(x)x^{-s-1}\,dx
\end{equation}
is holomorphic for $\Re s>1/(K+1)$.  Thus
\eqref{eq:D1-finite-expansion} gives a meromorphic continuation of
$\D_1$ to this half-plane.  Since $K$ is arbitrary, and the
continuations agree on their common domains by the identity theorem,
they define a meromorphic continuation of $\D_1$ to $\Re s>0$.

\Ecolor{In the half-plane $\Re s>1/(K+1)$, the only possible poles in
\eqref{eq:D1-finite-expansion} are $1/2,\ldots,1/K$. At $s=1/k$,
only the $k$th summand is singular, and its residue is $(k!)^{1/k}/k$.
This is nonzero, so each pole is simple. \Ecolor{Considering the right half-planes $\Re(s)>\frac{1}{k}$ for each $k$ one deduces that these are the only poles of $\D_1$ in the right half-plane.} 
As the poles of $\D_1$ accumulate at $0$, no
meromorphic continuation to a neighborhood of $0$ is possible.}
\end{proof}

\begin{remark}
The first residues are
\[
 \Res_{s=1/2}\D_1(s)=\frac1{\sqrt2},
 \qquad
 \Res_{s=1/3}\D_1(s)=\frac{6^{1/3}}3,
 \qquad
 \Res_{s=1/4}\D_1(s)=\frac{24^{1/4}}4.
\]
Thus each fixed nontrivial diagonal contributes one pole, and the
poles accumulate at the boundary point $s=0$.
\end{remark}

\section{Concluding remarks}

The results display a sharp change between the first and the higher exceptional sets.  Membership in $\E_2$ requires a collision between two distinct nontrivial diagonals, and fixed-pair finiteness reduces the global count below every positive power of $x$.  This tight density estimate also controls the complete overcount caused by all multiple intersections among any fixed finite collection of diagonals.  Consequently, for every fixed $K$, the first $K-1$ nontrivial diagonals give the expansion
\[
 A_1(x)=\sum_{k=2}^{K}(k!x)^{1/k}
 +O_K\!\left(x^{1/(K+1)}\right).
\]

In contrast, for every $H\geq2$, our estimates imply absolute and locally uniform convergence of $\D_H(s)$ throughout $\Re s>0$.

The bounds on the higher exceptional sets leave a natural quantitative gap. 

Our bound proves that for all $x \ge 3$, $A_2(x)<\exp\left(\frac{2\log x}{\left(\log\log x\right)^\frac{1}{4}}\right)$, while the only known repeated coefficients in $A_2(x)$ have $O(\log(\log(x)))$ density.
\Ecolor{  Finally, we believe that the bounds on $A_H(x)$ can be strengthened by combining Remond's techniques with specific properties of binomial Diophantine equations. }

\section*{Acknowledgments}
We thank the referee for carefully reading the paper and making a list of corrections. In particular, his suggestion that we use \cite{Remond2010} instead of \cite{BST} led to the effective bound in the second item of the main theorem. 

\section*
{AI disclosure}
The authors used OpenAI's ChatGPT 5.6 and Claude Opus 5.5 as interactive research and writing assistants during the development of this work, including testing and refining arguments, organizing proofs, and improving
the exposition. All mathematical
claims, proofs, references, and final formulations were independently reviewed
and verified by the authors, who take full responsibility for the content of
the paper.

\appendix
\section{A lower bound on $A_2(x)$}
\Ecolor{The Lind--Singmaster Fibonacci family gives a lower bound of order
$\log\log x$ for $A_2(x)$. We determine its density up to
a bounded error.}
Let $F_0=0$, $F_1=1$ and $F_{m+1}=F_m+F_{m-1}$.  Lind and Singmaster proved that, for every $i\geq1$,
\begin{equation}\label{eq:Fibonacci-identity}
 \binom{F_{2i+2}F_{2i+3}}{F_{2i}F_{2i+3}}
 =
 \binom{F_{2i+2}F_{2i+3}-1}{F_{2i}F_{2i+3}+1}.
\end{equation}
Denote the common value by $b_i$, and let
\begin{equation}\label{eq:F-family}
 \F=\{b_i:i\geq1\}.
\end{equation}
The two representations in \eqref{eq:Fibonacci-identity}, their symmetric partners, and the two edge representations show that $N(b_i)\geq6$.  Thus $\F\subseteq\E_2$.  The first value is $b_1=3003$, which has the additional representation $\binom{78}{2}$ and hence multiplicity eight.

Although $\F$ is infinite, it is exceptionally sparse.

\begin{proposition}\label{prop:F-density}
\Ecolor{Recall $\varphi=(1+\sqrt5)/2$.}  Then
\begin{equation}\label{eq:F-density}
 \#(\F\cap[1,x])
 =\frac{\log\log x}{4\log\varphi}+O(1).
\end{equation}
Consequently, for every $s>0$,
\begin{equation}\label{eq:F-series}
 \sum_{b\in\F}b^{-s}<\infty.
\end{equation}
\end{proposition}

\begin{proof}
Set
\[
 n_i=F_{2i+2}F_{2i+3},
 \qquad
 k_i=F_{2i}F_{2i+3}.
\]
Binet's formula for the Fibonacci sequence gives
\begin{equation}\label{eq:F-asymptotics}
 n_i=\frac{\varphi^{4i+5}}5(1+o(1)),
 \qquad
 \frac{k_i}{n_i}\longrightarrow\varphi^{-2}.
\end{equation}
Define the binary entropy function
\[
 h(t)=-t\log t-(1-t)\log(1-t)
\]
Stirling's formula, uniformly when $k_i/n_i$ stays in a compact subinterval of $(0,1)$, yields
\begin{equation}\label{eq:F-Stirling}
 \log b_i
 =n_i h(k_i/n_i)+O(\log n_i)
 =c\varphi^{4i}(1+o(1))
\end{equation}
for the positive constant
\[
 c=\frac{\varphi^5}{5}h(\varphi^{-2}).
\]
Therefore
\begin{equation}\label{eq:loglog-bi}
 \log\log b_i=4i\log\varphi+\log(c)+o(1).
\end{equation}
In particular, $b_i$ is strictly increasing for all sufficiently large $i$.  Solving \eqref{eq:loglog-bi} for $i$ proves \eqref{eq:F-density}.  The convergence in \eqref{eq:F-series} follows immediately from the double-exponential growth expressed by \eqref{eq:F-Stirling}.
\end{proof}

\bibliographystyle{amsplain}

\end{document}